\documentclass[onecolumn]{autart}    
\let\theoremstyle\relax

\usepackage[utf8]{inputenc}
\usepackage{comment}
\usepackage{algorithm}
\usepackage{algpseudocode}
\usepackage{breqn}
\usepackage{graphicx}
\usepackage[version=4]{mhchem}
\usepackage{siunitx}
\usepackage{longtable,tabularx}
\usepackage{graphics} % for pdf, bitmapped graphics files
\usepackage{float}
\usepackage{epsfig} 
\usepackage{times} 
\usepackage{amsmath} 
\usepackage{amssymb}  
\usepackage{mathtools}
\usepackage{pifont}
\usepackage[dvipsnames]{xcolor}
\usepackage{amsthm}
\usepackage{enumitem}
\usepackage{bm}
\usepackage{bbm}
\usepackage{soul}

\newcommand{\cov}{\text{cov}} % bob unbolded and lower-cased this

\newtheorem{assumption}{Assumption}

\newtheorem{lemma}{Lemma}

\newtheorem{theorem}{Theorem}

\DeclareMathOperator*{\argmin}{arg\,min}

\definecolor{paleGreen}{rgb}{.3, .7, .3}
\definecolor{coolBlue}{rgb}{.3, .5, 1}
\definecolor{rosePink}{rgb}{.9, .5, .4}
\definecolor{ghost}{rgb}{.8, .8, .8}

\newcommand{\ff}{\mbox{\textit{ff}}}
\newcommand{\mcbA}{\boldsymbol{\mathcal{A}}}
\newcommand{\mcbB}{\boldsymbol{\mathcal{B}}}
\newcommand{\mcbC}{\boldsymbol{\mathcal{C}}}
\newcommand{\mcbD}{\boldsymbol{\mathcal{D}}}

\begin{document}

\begin{frontmatter}
\title{Stable state and signal estimation in a network context} 
\author[UCSD]{Robert R. Bitmead}\ead{rbitmead@ucsd.edu}  
\address[UCSD]{Department of Mechanical \&\ Aerospace Engineering, University of California, San Diego, La Jolla CA 92093-0411, USA.}   
\thanks{This work was performed while the author was SimTech Visiting Professor at the Institute of Systems Theory, University of Stuttgart, partially supported by the Deutsche Forschungsgemeinschaft (DFG, German Research Foundation) under Germany's Excellence Strategy - EXC 2120/1 - 390831618 - EXC 2075/1 - 390740016.}

\begin{abstract}
Power grid, communications, computer and product reticulation networks are frequently layered or subdivided by design. The OSI seven-layer computer network model and the electrical grid division into generation, transmission, distribution and associated markets are cases in point. The layering divides responsibilities and can be driven by operational, commercial, regulatory and privacy concerns. From a control context, a layer, or part of a layer, in a network isolates the authority to manage, i.e. control, a dynamic system with connections into unknown parts of the network. The topology of these connections is fully prescribed but the interconnecting signals, currents in the case of power grids and bandwidths in communications, are largely unavailable, through lack of sensing and even prohibition. Accordingly, one is driven to simultaneous input and state estimation methods. This is the province of this paper, guided by the structure of these network problems. We study a class of algorithms for this joint task, which has the unfortunate issue of inverting a subsystem, which if it has unstable transmission zeros leads to an unstable and unimplementable estimator. Two modifications to the algorithm to ameliorate this problem were recently proposed involving replacing the troublesome subsystem with its outer factor from its inner-outer factorization or using a high-variance white signal model for the unknown inputs. The outer factor has only stable transmission zeros and so is stably invertible. Here, we establish the connections between the original estimation problem for state and input signal and the outputs/estimates from the algorithm applied solely to the outer factor. It is demonstrated that the state of the outer factor and that of the original system asymptotically coincide and that the estimate of the input signal to the outer factor has asymptotically stationary second-order statistics which are in one-to-one correspondence with those of the input signal to the original system, when this signal is itself stationary. Thus, the simultaneous input and state estimation algorithm applied just to the outer factor yields an unbiased state estimate for control and the statistics of the interface signals. We also show that the outer factor algorithm is the limit of the high-variance strategy, which yields an even simpler approach and implementation.
\end{abstract}
\end{frontmatter}

\section{Introduction}
\begin{figure}[ht]
\begin{centering}
\includegraphics[width=70mm]{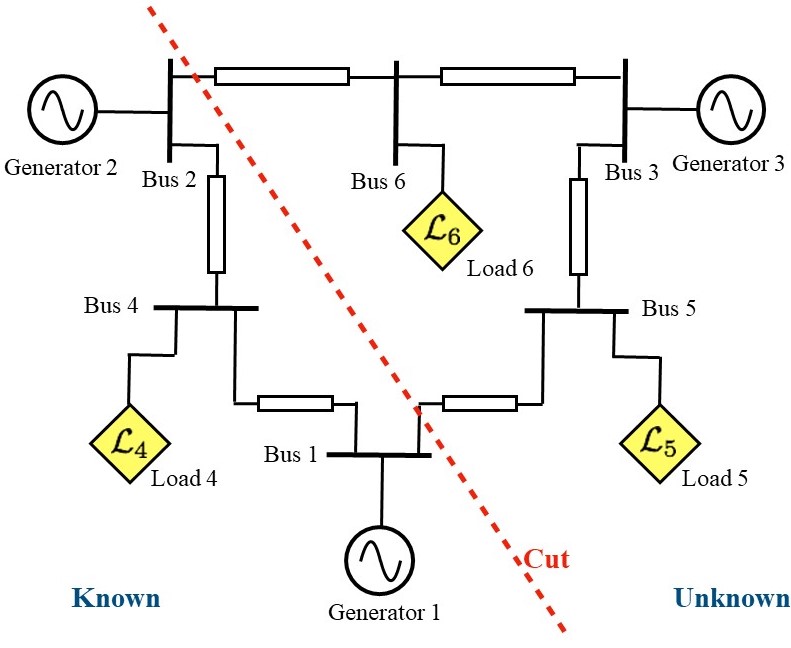}
\caption{ From \cite{AliMortenBobCCTA2019}: nine-bus power system with division into known part, Buses~1, 2 and 4 plus associated machines and Load~4, and unknown part by a circuit cut. Currents entering Buses~1 and 2 from the transmission lines on the unknown side comprise unknown input signal $d_t$.\label{fig:cutter}}
\end{centering}
\end{figure}

By design, networks involve multiple agents and operators who interact via the network links. Frequently, these multiple players are constrained to operate in only part of the overall network with responsibility solely for their section. Interactions with other denizens of the network occur through connections and the signals impinging there on the local section from these other operators. We assume that the operator knows their own system dynamics and the topology of interconnections to other operators. But they have no knowledge of the dynamics on the other side of these connections nor do they measure all the interconnection signals. In a power grid, where layering into generation, transmission, distribution, etc occurs, the generator operator should: know the dynamics of their generators and local synchronous machines; possess a set of available measurements from within their section; and,  have knowledge of where connections to other parts of the network occur, i.e. the buses where current is delivered to other levels of the network. Knowledge of these currents, unless they are directly measured, is not available nor is knowledge of the dynamics causing the currents outside the immediate section. We refer to this as the known part of the network. Such a network division is depicted for a simple power system in Figure~\ref{fig:cutter}.

In this scenario, the control of dynamic elements in the known part requires estimating their states and this is affected by the presence of disturbances arriving at the specified points of interconnection to the rest of the network. This joint problem of unbiased and least-squares estimation of both state and unknown input signal vectors has been studied under the general heading of Simultaneous Input and State Estimation (SISE) \cite{KitanidisAutom1987,GillijnsDeMoorAutom2007a,GillijnsDeMoorAutom2007b,SundaramHadjicostasAutom2008,YongZhuFrazzoliAutom2016,KongSukkariehAutom2019,gakisSmithIJC23}.

Our model of the known part of the network is
\begin{align}
P(z):\;\,x_{t+1}&=Ax_t+Gd_t+w_t,\,x_0,\label{eq:pxdef}\\
y_t&=Cx_t+Hd_t+v_t.\label{eq:pydef}
\end{align}
Since the system is linear, we have suppressed the presence of control and other known additive signals entering the system. Similarly, we have not included other modeled additive disturbances beyond $d_t$. Both can be directly incorporated using linearity.

In \eqref{eq:pxdef}-\eqref{eq:pydef}, $x_t\in\mathbb R^n$ is the state of our piece of the network, $d_t\in\mathbb R^m$ is the unknown input signal impinging on our part of the network from other sources, $w_t\in\mathcal R^n$ is the process noise taken to be white, zero mean, and of covariance $\mathcal Q$, $v_t\in\mathbb R^p$ is the measurement noise also taken to be white, zero-mean and of covariance $\mathcal R$, $y_t\in\mathbb R^p$ is the vector of measured system signals in the known part. Quantities $x_o$, $w_t$, $v_t$ and $d_t$ are independent. We take the system to be linear and time-invariant for clarity and because these requirements are needed for the stability properties of the SISE algorithms. 

There are variations on the SISE algorithm depending on the delay properties of $P(z)$. When $H=0$ and rank$(CG)=m$, the algorithm from \cite{GillijnsDeMoorAutom2007a} applies.
\begin{align}
X_t&=AP_{t-1}A^T+Q,\label{eq:xsise}\\
K_t&=X_tC^T(CX_tC^T+R)^{-1},\label{eq:sisek}\\
M_t&=[G^TC^T(CX_tC^T+R)^{-1}CG]^{-1}%\nonumber\\
%&\hskip 15mm\times 
G^TC^T(CX_tC^T+R)^{-1},\label{eq:sisem}\\
P_t&=(I-K_tC)\left[(I-GM_tC)X_t%\right.\nonumber\\
%&\hskip 8mm\left.\times
(I-GM_tC)^T+GM_tRM_t^TG^T\right]%\nonumber\\
%&\hskip 18mm
+K_tRM_t^TG^T,\label{eq:sisep}\\
\hat d_{t-1|t}&=M_t(y_t-CA\hat x_{t-1|t-1}),\label{eq:sised}\\
\hat x_{t|t}&=A\hat x_{t-1|t-1}+G\hat d_{t-1|t}+K_t%\nonumber\\
%&\hskip 8mm\times
(y_t-CA\hat x_{t-1|t-1}-CG\hat d_{t-1|t}),\label{eq:sisex}\\
\cov(x_t\vert&\mathbf Y^t)=P_t,\label{eq:sisecovx}
\end{align}
When rank$(H)=m$, a different variant is used. In more general cases, the ULISE algorithm of \cite{YongZhuFrazzoliAutom2016} or that of \cite{KongSukkariehAutom2019} based on strong detectability conditions can be used.

The important feature of the set of SISE algorithms, however, is that, since no explicit dynamic model is assumed for $d_t$, the methods proceed firstly by inverting $P(z)$ for $\hat d_{t-1|t}$ or $\hat d_{t|t}$ and thereafter estimating or simulating the state $x_t$ using the $d_t$ estimate and, if $p>m,$ the additional outputs. Accordingly, the stability of SISE algorithms is compromised when $P(z)$ possesses transmission zeros outside the unit circle \cite{YongZhuFrazzoliAutom2016,AbooshahabAlyaseenBitmeadHovd2022}. This instability can be ameliorated for these algorithms only by having sufficient additional measurements, i.e. $p$ being sufficiently greater than $m$, that a detectability condition holds.

However, in \cite{AbooshahabAlyaseenBitmeadHovd2022} two workarounds are proposed to guarantee stability of a modified SISE even when $P(z)$ has transmission zeros outside the unit circle.
\begin{description}
\item[High-$\mathcal D$ Kalman filtering:] It is shown in \cite{BitmeadHovdAbooshahab:Autom2019} that SISE coincides precisely with the Kalman filter with $\mathcal R>0$ and $d_t$ taken as white, independent from $x_0$, $w_t$ and $v_t$ and with covariance $\mathcal D$ with $\mathcal D^{-1}=0$. That is, SISE is singular Kalman filtering \cite{ShakedTAC1985,PrielShakedCDC1986}. The proposal in \cite{AbooshahabAlyaseenBitmeadHovd2022} is to implement the standard Kalman filter based on white $d_t$ with large $\mathcal D$. That is, with process noise covariance $G\mathcal DG^\top+\mathcal Q$. This filter is guaranteed stable under the usual conditions of $\mathcal R>0$ and $[A,\mathcal Q]$ stabilizable. Although this comes at the price of small (dependent on $\mathcal D^{-1}$) bias in the estimation of $d_t$.
\item[Inner-outer factorization:]  Plant $P(z)$ is factored into its inner and outer parts $P(z)=P_o(z)P_i(z),$ with $P_i(z)$ being inner, i.e. a stable all-pass function, and $P_o(z)$ being outer, i.e. stable with stable transmission zeros. The SISE algorithm \eqref{eq:xsise}-\eqref{eq:sisecovx} is applied to $P_o(z)$ in place of $P(z)$. The signal $f_t=P_i(s)d_t$ is estimated in place of $d_t$.
\end{description}
The contribution of this paper is to show that:
\begin{enumerate}[label=(\roman*)]
\item the states of $P(z)$ and $P_o(z)$ asymptotically coincide so the state estimate, $\hat x_{t|t}$, from stable SISE applied to $P_o(z)$ is an asymptotically unbiased estimate of the state of $P(z)$;
\item when $d_t$ is stationary, the second-order statistics of $d_t$ are simply recoverable from those of the input $f_t$ estimates from $P_o(z)$;
\item the two workaround methods coincide.
\end{enumerate}

\section*{Notation}
%\textit{Discrete time, complex variable $z$.}
RH$^\infty$ denotes the set of proper rational matrices analytic in $|z|\geq 1.$ For rational matrix $P(z),$ its paraconjugate is $P(z)^\sim=P\left(z^{-1}\right)^*$, where $\cdot^*$ denotes Hermitian conjugate. Denote the transpose operation by $P^\top(z)$. For $P(z)$ in RH$^\infty$, $P(z)$ is inner if $P(z)P(z)^\sim=I$ and $P(z)$ is outer if it has full column rank for every $z$ in $|z|\geq1$. The MacMillan degree, $\delta(P)$, of transfer function $P(z)$ is the minimal state dimension for a realization of $P(z)$ \cite{NewcombBook1966,macfarlaneKarcaniasIJC1976,Kailath:80}. Transfer function $P(z)$ is \textit{regular} if $\delta(PP^\sim)=2\delta(P)$. Transfer function $R(z)$ is a \textit{spectral factor} of $PP^\sim$ if $RR^\sim=PP^\sim$ with $R$ in RH$^\infty$ and column rank\,$R(z)$ full for all $z$ in $|z|\geq 1.$ Note that, since we are focused on estimation, our definitions of inner, outer, spectral factorization and inner-outer factorization are transposed from their usual definitions in robust control.

\section{Inner-outer factorization for state estimation}
For $p\times m$ transfer function $P(z)$ in RH$^\infty$, inner-outer factorization (for state estimation) writes $P(z)=P_o(z)P_i(z)$ with $P_i(z)$ being $r\times m,\,r\leq p,$ and inner, i.e. $P_i(z)P_i(z)^\sim=I_r$, and $P_o(z)$ being $p\times r$ and outer. 
A scalar clarifying example is in order. Consider
\begin{align*}
P(z)&=\frac{(z-2)(z-3)(z-0.9)(z-0.8)}{(z-\half)(z-0.7)(z+j\half)(z-j\half)}=\overbrace{\frac{6(z-\frac{1}{3})(z-0.9)(z-0.8)}{(z-0.7)(z+j\half)(z-j\half)}}^{P_o(z)}\times\overbrace{\frac{(z-2)}{2(z-\half)}\times\frac{(z-3)}{3(z-\frac{1}{3})}}^{P_i(z)}
,\\
 P_i(z)P_i(z)^\sim&=\frac{1}{6}\frac{(z-2)(z-3)}{(z-\half)(z-\frac{1}{3})}\frac{1}{6}\frac{(z^{-1}-2)(z^{-1}-3)}{(z^{-1}-\half)(z^{-1}-\frac{1}{3})}
=\frac{1}{36}\frac{(z-2)(z-3)}{(z-\half)(z-\frac{1}{3})}\frac{6(z-\frac{1}{3})(z-\half)}{\frac{1}{6}(z-3)(z-2)}=1.
\end{align*}
The outer factor, $P_o$, has only stable poles and zeros and has $\delta(P_o)=3$. The inner factor has stable poles and consists of two terms. The first, $\frac{z-2}{2(z-\half)},$ is already an inner factor of $P$ -- a \textit{free inner factor} according to Green \cite{GreenInnerOuterSCL1988}. Its presence is indicated because $\delta(PP^\sim)=6<2\delta(P)=8.$ That is, this $P$ is not \textit{regular} as defined above. The inner term, $\frac{z-3}{3(z-\frac{1}{3})},$ is a non-free inner factor of $P$.

With the ordering $P=P_iP_o$, this is a standard calculation in robust control \cite{GreenInnerOuterSCL1988,ChenFrancisSIAMMAA1989,Green&Limebeer:95,Zhou&Doyle&Glover:95}. Our ordering is simply computed by applying the standard algorithm to $P^\top(z)$ and the transposing the factorization.

We make the following assumption
\begin{assumption}\label{ass:sys}
\begin{enumerate}[label=(\roman*)]
\item $P(z)$ is stable, i.e. $P(z)\in\text{RH}^\infty$.\label{pt:stable}
\item $P(z)$ is regular, i..e. $\delta(PP^\sim)=2\delta(P)$. 
\item Realization $[A,G,C,H]$ is minimal.
\item $[A,\mathcal Q]$ is reachable, $\mathcal R>0$.
\end{enumerate}
\end{assumption}
%We shall later discuss relaxing the stability condition Assumption~\ref{ass:sys}.\ref{pt:stable}.

We have the following from the Appendix, converted to discrete time and transposed to reorder the factors and, thereby, to focus on observability in place of reachability. Here $\ell$ is the number of unstable transmission zeros of $P$.
\begin{theorem}[From Appendix]\label{thm:cancelo}
For regular $P(z)$ in RH$^\infty$ and inner-outer factorization $P(z)=P_o(z)P_i(z)$ with $\delta(P_i)=\ell$,
\begin{enumerate}[label=(\roman*)]
\item If $P(z)$ has minimal state-variable realization $H+C(zI-A)^{-1}G$ then $P_o(z)$ has minimal state-variable realization $H_o+C(zI-A)^{-1}G_o$. That is, with identical $A$ and $C$ to those of  $P$.\label{pt:minoto}
\item The state-space realization of the inner factor $P_i(z)=U^\top+\hat B^\top(sI-\hat A^\top)\hat C^\top$ is minimal.\label{pt:minino}
\item $n=\delta(P)=\delta(P_o)\leq \delta(P_o)+\delta(P_i)=n+\ell,$ with equality only when $\delta(P_i)=0$, i.e. $P=P_o$.\label{pt:delo}
\item The full state of the outer factor $P_o$ is observable\footnote{Note that this usage of observability for stochastic systems relies on the formal definition in \cite{LiuBitmead_Autom2011}, which subsumes the more usual deterministic concept of observability. Particularly, the role of the system input signal in state estimation is central to understanding observability.} from output $y_t$ and an $\ell$-dimensional subspace, comprising the states of the inner factor $P_i(z),$ is unobservable from the output of $P$.\label{pt:unreo}
\item The unobservable modes are at the eigenvalues of $\hat A^\top$, which are stable by construction and lie at the inverses of the transmission zeros of $P(s)$ outside the unit disc.\label{pt:stblo}
\end{enumerate}
\end{theorem}

\section{State estimation using the outer factor}
Consider the plant given by \eqref{eq:pxdef}-\eqref{eq:pydef}
\begin{align*}
P(z):\;\,x_{t+1}&=Ax_t+Gd_t+w_t,\,x_0,\\
y_t&=Cx_t+Hd_t+v_t,
\end{align*}
and its outer-factor
\begin{align}
P_o(z):\;\,x^o_{t+1}&=Ax^o_t+G_of_t+w_t,\,x^o_0,\label{eq:poxdef}\\
y_t&=Cx^o_t+H_of_t+v_t,\label{eq:poydef}
\end{align}
where $f_t$ is the output of the inner factor driven by $d_t$.
\begin{align}
P_i(z):\;\;x^i_{t+1}&=A_ix^i_t+G_id_t,\;\;x^i_0,\label{eq:pixdef}\\
f_t&=C_ix^i_t+H_id_t.\label{eq:piydef}
\end{align}

\begin{lemma}\label{lem:stests}
Denote the observability Gramian of $[A,C]$ as
\begin{align}\label{eq:obgramdef}
W_o(N)&=\sum_{j=0}^{N-1}{A^{\top^j}C^\top CA^j}.
\end{align}
Define the signals
\begin{align}
\acute y_t&=y_t-\left (Hd_t+v_t+\sum_{j=0}^{N-1}{CA^{N-1-j}(Gd_{t-j}+w_{t-j})}\right),\label{eq:acutedef}\\
\grave y_t&=y_t-\left (H_of_t+v_t+\sum_{j=0}^{N-1}{CA^{N-1-j}(G_of_{t-j}+w_{t-j})}\right).\label{eq:gravedef}
\end{align}
Then, for any $t\geq N\geq n$, the state dimension, $W_o(N)$ is invertible and
\begin{align}
x_{t-N+1}&=W_o(N)^{-1}\sum_{j=0}^{N-1}{A^{\top^{N-1-j}}C^\top \acute y_{t-j}}.,\label{eq:statin}\\
x^o_{t-N+1}&=W_o(N)^{-1}\sum_{j=0}^{N-1}{A^{\top^{N-1-j}}C^\top \grave y_{t-j}}.\label{eq:ostatfin}
\end{align}
\end{lemma}
\begin{proof} The invertibility of the observability Gramian follows from the minimality of the state-variable realizations for $P$ and $P_o$.
From the state equation, we have
\begin{align*}
\begin{bmatrix}\acute y_{t-N+1}\\\acute y_{t-N+2}\\\vdots\\\acute y_{t}\end{bmatrix}&=\begin{bmatrix}C\\CA\\\vdots\\CA^{N-1}\end{bmatrix}x_{t-N+1},\\
\begin{bmatrix}C^\top&A^\top C^\top&\dots&A^{\top^{N-1}}\end{bmatrix}\begin{bmatrix}\acute y_{t-N+1}\\\acute y_{t-N+2}\\\vdots\\\acute y_{t}\end{bmatrix}&=
\begin{bmatrix}C^\top&A^\top C^\top&\dots&A^{\top^{N-1}}\end{bmatrix}\begin{bmatrix}C\\CA\\\vdots\\CA^{N-1}\end{bmatrix}x_{t-N+1},\\
\sum_{j=0}^{N-1}{A^{\top^{N-1-j}}C^\top \acute y_{t-j}}.&=W_o(N)x_{t-N+1}.
\end{align*}
The proof for the $x^0_t$ calculation is analogous.
\end{proof}

\begin{lemma}\label{lem:same}
The inner-outer factorization $P(s)=P_o(s)P_i(s)$ and unobservability of $x^i_{t}$ imply that
\begin{align*}
Hd_t+\sum_{j=0}^{N-1}{CA^{N-1-j}Gd_{t-j}}=H_of_t+\sum_{j=0}^{N-1}{CA^{N-1-j}G_of_{t-j}}+\kappa_t\alpha^t,\\
\end{align*}
where $|\kappa_t|$ is uniformly bounded and $\alpha$ is the maximal modulus eigenvalue of $\hat A$.
\end{lemma}
\begin{proof}
The output, $y_t,$ of $P(z)$ with input signal $d_t$ is
\begin{align*}
y_t&=CA^tx_0+Hd_t+\sum_{j=0}^{t-1}{CA^{t-j}Gd_j}.
\end{align*}
Similarly, the output, $y^o_t,$ of $P_o(z)$ with input $f_t=C_iA_i^tx^i_0+H_id_t+\sum_{j=0}^{t-1}{C_iA_i^{t-j}G_id_j}$ is
\begin{align*}
y^o_t&=CA^tx^o_0+H_of_t+\sum_{j=0}^{t-1}{CA^{t-j}G_of_j},\\
&=CA^tx^o_0+Hd_t+\sum_{j=0}^{t-1}{CA^{t-j}Gd_j}.
\end{align*}
This latter equality follows two reasons.
\begin{enumerate}
\item The factorization, $P(z)=P_o(z)P_i(z),$ states that the convolution of the impulse response of $P_o$ with that of $P_i$ yields the impulse response of $P$, term by term. So, the zero-state response of $P_o$ driven by the zero-state response of $P_i$ with input $d_t$ is identical with the zero-state response of $P$ driven by $d_t$. What remains to be handled are the responses to the non-zero initial states, $x_0$, $x^o_0$ and $x^i_0$.
\item The construction in the Appendix to yield \eqref{eq:omgstate}-\eqref{eq:omgout} shows that the impulse response of the transformed non-minimal state realization \eqref{eq:omgstate}-\eqref{eq:omgout} has $\ell$-dimensional subspace $x^i_t$ unobservable. This state component does not affect $y_t$ and evolves according to $\hat A^\top$. \end{enumerate}
\end{proof}

\begin{theorem}\label{thm:same}
As $t\to\infty$, the state, $x^o_t$, of $P_o(z)$ converges exponentially fast to the state, $x_t,$ of $P(z)$ with a rate determined by the maximal modulus eigenvalue of $A$.
\end{theorem}

\begin{proof}
From Lemma~\ref{lem:same} and using \eqref{eq:acutedef}-\eqref{eq:gravedef}, the signals $\acute y_t$ and $\grave y_t$ converge exponentially fast in $t$. Thus, the states, which satisfy identities \eqref{eq:statin} and \eqref{eq:ostatfin}, converge.
\end{proof}

The upshot of this result is that the state of $P(z)$ driven by $(d_t,w_t)$ and the state of $P_o(z)$ driven by $(f_t,w_t)$ asymptotically coincide. Hence, using SISE on $y_t$ to produce estimates of $f_t$ and $x^o_t$ yields estimates of $f_t$ and of the state $x_t$.

\section{Input estimation using the outer factor}
The results of the earlier section show that the states of $P(z)$ and $P_o(z)$ are asymptotically identical. So, applying the SISE estimator to $P_o(z)$, since it has no unstable transmission zeros, yields a stable estimator for the state, $x_t$, of $P(z)$. The input estimated by this well behaved SISE filter is, however, not $d_t$ but its filtered version $f_t=P_i(z)d_t.$ Since $P_i$ is inner, there are direct invertible relations between the second-order statistics of $f_t$ and those of $d_t$, as will be derived shortly. Although, since $P_i$ is not stably invertible, one cannot recover $d_t$ from $f_t$ by stable filtering.

As we outline in the introduction, the network context can dictate that direct measurement of connection signals cannot be made for privacy or commercial concerns. However, the statistics of these signals are required for plant sizing and operational reasons. The following theorem demonstrates the connection between the second-order statistics of $f_t$ and those of $d_t$.

Standard results on stationary stochastic processes \cite{GrayDavission2004} yield the following connections.
\begin{theorem}
 Suppose stationary $m$-vector signal $d_t$ has mean value, autocovariance function and power spectral density as follows.
 \begin{align*}
 \bar d&={\rm E}(d_t),\\
 R_{dd}(\tau)&={\rm E}\left[\left(d_t-\bar d\right)\left(d_{t+\tau}-\bar d\right)^\top\right],\\
 \Phi_{dd}(\omega)&=\mathcal F\left[R_{dd}(\tau)\right],
% &\triangleq \mathcal F^{-1}\left[\Phi_{dd}(s)\right].
 \end{align*}
 where $\mathcal F[\cdot]$ is the discrete-time Fourier transform. For inner filter $P_i(z)$, the $r$-vector filtered signal $f_t=P_i(z)d_t$ is asymptotically stationary with mean, power spectral density and covariance
\begin{align*}
\bar f\triangleq {\rm E}(f_t)
 &=P_i(1)\bar d,\\
 \Phi_\ff(\omega)&=\mathcal F\left[{\rm E}\left[\left(f_t-\bar f\right)\left(f_{t+\tau}-\bar f\right)^\top\right]\right],\\
 &=P_i(e^{j\omega})\Phi_{dd}(\omega)P_i^\top(e^{-j\omega}),\\
 R_\ff(\tau)&=\mathcal F^{-1}\left[\Phi_\ff(\omega)\right].
  \end{align*}
 Further,
 \begin{align*}
 \bar d&=P_i^\top(1)\bar f,\\
% &=P_i^\top(0)\bar d,\;\text{and}\\
 \Phi_{dd}(\omega)&=P^\top_i(e^{-j\omega})\Phi_\ff(\omega)P_i(e^{j\omega}),\\
 R_{dd}(\tau)&=\mathcal F^{-1}\left[\Phi_{dd}(\omega)\right].
  \end{align*}
 So, the steady-state second-order statistics of $d_t$ are simply recovered from those of $f_t$ and vice versa.
 \end{theorem}

\section{High-$\mathcal D$ filtering and outer factors}\label{sec:inest}
The high-$\mathcal D$ filter with $p$, the number of outputs, greater than $m$, the number of $d_t$ inputs, may be solved using the sequential decomposition of Priel and Shaked for the \textit{partially singular} filtering problem \cite{PrielShakedCDC1986}. The decomposition transforms the output into two signals $\begin{bmatrix}y_1^\top&y_2^\top\end{bmatrix}^\top=\begin{bmatrix}C_1^\top&C_2^\top\end{bmatrix}^\top x$ with $m\times m$ matrix $C_1G$ invertible. The process noise variance for the filtering problem with output signal $y_1$ is $G\mathcal D G^\top+\mathcal R$. When $\mathcal D$ increases without bound, the $\mathcal R$ term becomes unimportant and we are left with a fully singular filtering problem, whose solution is known. The other measurements, $y_2$, can then be incorporated in a standard fashion, effectively as a further measurement update to the state estimate based on $y_1$.

To appreciate the invertible-$C_1G$ solution, it is useful to consider the Return Difference Equality or spectral factorization formulation of optimal estimation. The Kalman filtering version of this equality is derived in \cite[Chapter~5]{bi:BGW90}. For a plant described by matrices $A$ and $C$ with process noise covariance $R$ and measurement noise covariance $R$,
\begin{align}\label{eq:edr}
R +C(zI-A)^{-1}Q(z^{-1}I-A^\top)^{-1}C^\top&=[I+C(zI-A)^{-1}L](C\Sigma C^\top+R)[I+L^\top(z^{-1}I-A^\top)^{-1}C^\top],
\end{align}
where $\Sigma$ is the positive definite solution of the prediction Algebraic Riccati Equation, i.e. the covariance of the state prediction error, and $L$ is the Kalman predictor gain $A\Sigma C^\top(C\Sigma C^\top +R)^{-1}.$

Note also the identity regarding the right-hand-side factor above.
\begin{align}
\left[I+C(zI-A)^{-1}L\right]^{-1}&=I-C(zI-A+LC)^{-1}L.\label{eq:facto}
\end{align}
This latter transfer function is that from output measurement $y_t$ to innovations sequence in the Kalman predictor.
\begin{align*}
\hat x_{t|t-1}&=(A-LC)\hat x_{t|t-1}+Ly_t,\\
y_t-\hat y_{t|t-1}&=y_t-C\hat x_{t|t-1}.
\end{align*}
So, the Return Difference Inequality \eqref{eq:edr} links the additive terms on the left side dealing with the estimation problem statement, $[A,C,Q,R],$ and the multiplicative terms on the right-hand side, which define gain $L$ and prediction error covariance $\Sigma$. The factoring of the additive left side to yield the multiplicative right side is the link between spectral factorization and optimal filtering \cite{Anderson&Moore:79}.

For $y_1=C_1x$ as the measurement $m$-vector and as pointed out in \cite{PrielShakedCDC1986}, the filtered error covariance, $S_t$, and the predicted error covariance, $\Sigma_t$, satisfy $S^{-1}=\Sigma^{-1}+C_1R^{-1}C_1^\top$ and $\Sigma=ASA^\top+G\mathcal D G^\top$. With $R$ fixed and $\mathcal D$ unbounded, $C_1\Sigma C_1$ tends to $C_1G\mathcal DG^\top C_1^\top$ and \eqref{eq:edr} becomes.
\begin{align*}
R +C_1(zI-A)^{-1}G\mathcal D G^\top(z^{-1}I-A^\top)^{-1}C_1^\top&=[I+C_1(zI-A)^{-1}L](C_1G\mathcal DG^\top C_1^\top+R)[I+L^\top(z^{-1}I-A^\top)^{-1}C^\top].
\end{align*}
Take $C_1(zI-A)^{-1}G=P(z)=P_o(z)P_i(z)$  and $\mathcal D=\epsilon^{-1}I_m$ for $\epsilon\to 0$ (In \cite{BitmeadHovdAbooshahab:Autom2019}, the high-$\mathcal D$ filter is shown to be identical for any $\mathcal D$ being rank $m$ and $\mathcal D^{-1}\to 0$.) then we have
\begin{align*}
P_o(z)P_i(z)\mathcal DP_i^\sim(z)P_o^\sim(z)&=[I+C_1(zI-A)^{-1}L]C_1G\mathcal DG^\top C_1^\top[I+L(zI-A)^{-1}C]^\sim,\\
\epsilon^{-1}P_o(z)P_i(z)P_i^\sim(z)P_o^\sim(z)&=\epsilon^{-1}[I+C_1(zI-A)^{-1}L]C_1GG^\top C_1^\top[I+L(zI-A)^{-1}C]^\sim,\\
P_o(z)P_o(z)^\sim&=\left\{[I+C_1(zI-A)^{-1}L]C_1G\right\}\left\{[I+C_1(zI-A)^{-1}L]C_1G\right\}^\sim.
\end{align*}
The corresponding stable-transmission-zero spectral factor corresponds to $P_o(z)$.
From this fully singular filter, one then proceeds per \cite{PrielShakedCDC1986} to develop the partially singular filter with additional measurement $y_2$.

\begin{theorem}\label{thm:singular}
Fixing $\mathcal R>0$ and taking $Q=G\mathcal DG^\top+\mathcal Q=\epsilon^{-1}GDG^T+\mathcal Q$ with $D$ rank$\,m$ and $\epsilon\to 0$ yields the Kalman filter for $P_o(z)$, the outer factor of $P(s)$. This filter is stable.
\end{theorem}

\textbf{Remark}
There are control counterparts, which are analyzed in revealing detail in \cite{Astrom&Wittenmark:97} for the minimum-variance scalar control problem, the dual to our singular filtering problem here. Other authors have written in detail about singular optimal control, see \cite{BellJacobson1975,ClementsAnderson1978}. In \cite{AliMortenBobCCTA2019} $\mathcal D$ is taken as $10^6\,I_m$ for a power systems example. Such a value did not create solution accuracy issues for the algebraic Riccati equation.

\section{Conclusions}
The principal objective of this paper is to tie together and provide theoretical support for the SISE algorithm workarounds to avoid stability problems in applying these methods in network contexts, where the information architecture requires simultaneous input and state estimation, i.e. SISE, in order to effect local control and/or estimate interactions.  Input estimation problems abound in many domains, notably in instrument deconvolution \cite{fang2013simultaneous,gakisSmithIJC23}.
But networks would appear to be especially fruitful because of the independence and heterogeneity of the agents operating across the networks and the absence of knowledge of their motives and behaviors. However, the known stability issue of SISE algorithms is an immediate impediment.

The results in this paper validate the use of guaranteed stable SISE estimators using the outer-factor of the plant system. The states of the original system and of its outer factor are proven to be asymptotically identical. Further, the all-pass-filtered interconnection signal, $f_t,$ constructed by this stable method is shown to maintain the second-order statistics of the original signals, $d_t$. So, this application might even have privacy advantages.

The final piece of the puzzle, that the inner-outer factorization and the high-$\mathcal D$ approaches coincide, actually suggests avoiding the problematic SISE algorithms all together and using the standard Kalman filter with a special choice of process noise covariance. These approaches have been applied successfully in partially known power systems \cite{AliMortenBobCCTA2019}.

\section*{Appendix -- On Green's inner-outer factorization, $P=P_iP_o,$ and its state estimation variant, $P=P_oP_i$}
We take the slightly circuitous route of transforming to continuous time for the analysis of inner-outer factorization. We do this because, despite the existence of many papers detailing the discrete-time calculation of inner-outer factorizations, see e.g. \cite{IonescuOaraTAC1996,IonescuOaraSICON1996,LinChenSaberiShamashTransCSI1996,ChuSpectralACC1988}, Green's work in continuous-time  focuses on delivering minimal realizations of each factor. We appeal to this property to analyze reachability here. Note, Green computes the usual (control) inner-outer factorization ordering $P(s)=P_i(s)P_o(s),$ which differs from ours. Also, in continuous time, RH$^\infty$ denotes the set of proper rational matrices analytic in Re$(s)>0.$ For rational matrix $P(s),$ $P^\sim(s)=P^*(-\bar s)$ is the paraconjugate with $\bar\cdot$ being complex conjugation. For $P(s)$ in RH$^\infty$, $P(s)$ is inner if $P^\sim(s)P(s)=I$ and $P(s)$ is outer if it has full row rank for every $s$ in Re$(s)>0$. Corresponding definitions of MacMillan degree, regular, and spectral factor apply. We rely on the Tustin transformation between the two time domains to carry these results back to discrete time.

\section*{A.I Tustin transformation between discrete and continuous time domains}
The Tustin transform, an instance of the more general M\"obius transformation, maps between the discrete-time complex $z$ plane and the continuous-time complex $s$-plane and is scaled by $\omega_0=\frac{1}{T}$ where $T$ is the sampling time.
\begin{align}\label{eq:tustin}
z=\frac{\omega_0+s}{\omega_0-s},\text{ and   }s=\omega_0\frac{z-1}{z+1}.
\end{align}
Proper discrete state-space transfer function $H(z)$ transforms to proper continuous transfer function $G(s)$ as follows.
\begin{align}\label{eq:tustmats}
H(z)=D+C(zI-A)^{-1}B \implies G(s)=\bar D+\bar C(sI-\bar A)^{-1}\bar B\;\;\;\text{with}\;\;
\begin{cases}
\bar A=\omega_0(A-I)(A+I)^{-1},\\
\bar B=\sqrt{2\omega_0}(I+A)^{-1}B,\\
\bar C=\sqrt{2\omega_0}C(I+A)^{-1},\\\bar D=D-C(I+A)^{-1}B.
\end{cases}
\end{align}
\begin{lemma}\label{lm:tustin}
The Tustin transformation preserves innerness and outerness.
\end{lemma}
\begin{proof} The Tustin transformation maps the inside of the unit disc in the $z$-plane to the left half-plane in $s$, and vice versa. So stability is preserved under the transformation, as is outerness, since these properties relate solely to the location of poles and zeros with relation to the stability boundary.

With $H(z)=D+C(zI-A)^{-1}B$,
\begin{align*}
H(z^{-1})&=(D-CA^{-1}B)-CA^{-1}(zI-A^{-1})^{-1}A^{-1}B.
\end{align*}
Applying \eqref{eq:tustmats}, the corresponding transformed $\breve G(s)=\breve D+\breve C(sI-\breve A)^{-1}\breve B$ will have
\begin{align*}
\breve A&=\omega_0(A^{-1}-I)(A^{-1}+I)^{-1}=-\omega_0(A-I)(A+I)^{-1},\\
\breve B&=\sqrt{2\omega_0}(I+A^{_1})^{-1}=\sqrt{2\omega_0}(I+A)^{-1}B,\\
\breve C&=-\sqrt{2\omega_0}CA^{-1}(I+A^{-1})^{-1}=-\sqrt{2\omega_0}C(I+A)^{-1},\\
\breve D&=D-CA^{-1}B+CA^{-1}(I+A^{-1})^{-1}A^{-1}B,\\
&=D-C\left[A^{-1}-A^{-1}(I+A^{-1})^{-1}A^{-1}\right]B,\\
&=D-C\left[A^{-1}-(I+A)^{-1}A^{-1}\right],\\
&=D-C(I+A)^{-1}\left[(I+A)A^{-1}-A^{-1}\right]B,\\
&=D-C(I+A)^{-1}B.
\end{align*}
That is, $\breve G(s)=G(-s)$. So, if $H(z)^\sim H(z)=I$ then $G(s)^\sim G(s)=I$. So innerness is also preserved.
\end{proof}
\section*{A.II Green's algorithm from robust control: $P(s)=P_i(s)P_o(s)$}
Michael Green \cite{GreenInnerOuterSCL1988} provides an explicit computational algorithm for the calculation of the inner-outer factorization $P(s)=P_i(s)P_o(s)$ with $P$ in RH$^\infty$. Inner-outer factorization is unique up to inclusion of a arbitrary unitary matrix between the terms. Green's algorithm begins with a minimal state-space realization of $P(s)$ and yields minimal state-space realizations of both $P_i(s)$ and $P_o(s)$. We note, however, that, in general, the composite state-space realization of the product $P_i(s)P_o(s)$ is non-minimal but stabilizable; a property we exploit in the paper. Green distinguishes between transfer functions with \textit{free inner factors}, as explained earlier in the paper, and regular transfer functions. Our problem of network estimation assumes that the transfer function $P$ is regular. Otherwise, there is a dimension mismatch between then state of $P_o$ and that of $P$.

Take $P(s)=D+C(sI-A)^{-1}B$ as a minimal state-space realization of $m\times p$ transfer function $P(s)$, which we assume to be stable and regular.
\begin{enumerate}
\item Compute the observability Gramian $Q=Q^*>0$ of $P$ satisfying
\begin{align*}
QA+A^*Q+C^*C&=0.
\end{align*}
\item Compute $P_o(s)$ as the spectral factor of $P^\sim P$. Since $P$ is regular, $\delta(P_o)=\delta(P)$. The dimensions of $P_o$ are $r\times p$ with $r\leq m$. This factor has the following minimal state-space realization.
\begin{align*}
P_o(s)&=J+H(sI-A)^{-1}B.
\end{align*}
\item Compute the observability Gramian $X=X^*>0$ of $P_o$ satisfying
\begin{align*}
XA+A^*X+H^*H=0.
\end{align*}
{}[Note that $X\leq Q$ and $\ell=$rank$(X-Q)$ is the number of transmission zeros of $P(s)$ in Re$(s)>0$.]
\item Find unitary (for us, orthogonal) transformation, $V$ with $V^*=V^{-1}$, such that
\begin{align*}
V(Q-X)V^*=\begin{bmatrix}0&0\\0&\Sigma\end{bmatrix},
\end{align*}
with $\ell\times \ell$ matrix $\Sigma>0$. Many matrix decompositions yield this.
\item Transform and partition the state matrices as follows.
\begin{align*}
VAV^{-1}&=\begin{bmatrix}A_{11}&A_{12}\\A_{21}&A_{22}\end{bmatrix},
VB=\begin{bmatrix}B_1\\B_2\end{bmatrix},
CV^{-1}=\begin{bmatrix}C_1&C_2\end{bmatrix}, HV^{-1}=\begin{bmatrix}H_1&H_2\end{bmatrix}.
\end{align*}
So that $A_{11}$, $B_1$ have $n-\ell$ rows, $C_1$ and $H_1$ have $n-\ell$ columns, etc.
\item Find $p\times r$ matrix $U$ such that $U^*U=I_r$ and 
\begin{align*}
\begin{bmatrix}C_1&D\end{bmatrix}=U\begin{bmatrix}H_1&J\end{bmatrix}.
\end{align*}
\item Then
\begin{align*}
\hat C&=UH_2-C_2&(p\times\ell),\\
\hat B&=\Sigma^{-1}(C_2^*U-H_2^*)&(\ell\times r),\\
\hat A&=A_{22}+\hat BH_2&(\ell\times\ell),
\end{align*}
yields $p\times r$ 
\begin{align}
P_i(s)&=U+\hat C(sI-\hat A)^{-1}\hat B&(p\times r),\label{eq:gidef}
\end{align}
with $P_i(s)$ inner and $P(s)=P_i(s)P_o(s)$.
\end{enumerate}

Using Green's construction, we have the following result.
\begin{theorem}[Green \cite{GreenInnerOuterSCL1988} plus minor extensions]\label{thm:cancel}
For regular $P(s)$ in RH$^\infty$ and inner-outer factorization $P(s)=P_i(s)P_o(s)$ as above with $\delta(P_i)=\ell$, the following properties hold.
\begin{enumerate}[label=(\roman*)]
\item If $P(s)$ has minimal state-variable realization $D+C(sI-A)^{-1}B$ then $P_o(s)$ has minimal state-variable realization $J+H(sI-A)^{-1}B$.\label{pt:minot}
\item The state-space realization \eqref{eq:gidef} of inner factor, $P_i(s)=U+\hat C(sI-\hat A)^{-1}\hat B,$ is minimal.\label{pt:minin}
\item $\delta(P)=\delta(P_o)\leq \delta(P_o)+\delta(P_i),$ with equality only when $\delta(P_i)=0$.\label{pt:del}
\item The full state, $x^o$, of $P_o(s)$ is reachable from the input to $P(s)$ and an $\ell$-dimensional subspace, comprising the sum of the states of $P_i$ and a linear combinations of the states of $P_o,$ is unreachable from the input to $P(s)$.\label{pt:unre}
%The states of the inner factor $P_i(s)$ are unreachable from the input to $P(s)$.\label{pt:unre}
\item The unreachable modes lie at the eigenvalues of $\hat A$ and so are stable. These mode values are the negatives of the right half-plane transmission zeros of $P(s)$. The reachable modes are at the eigenvalues of $A$.\label{pt:stb|}
\end{enumerate}
\end{theorem}
\begin{proof}
Green establishes Part~\ref{pt:minot} for regular $P(s)$. Part~\ref{pt:minin} is his central result. Part~\ref{pt:del} follows from the regularity of $P(s)$ and the properties of MacMillan degree. For Part~\ref{pt:unre}, consider the state-variable realization of $P_i(s)P_o(s)$.
\begin{align*}
\begin{bmatrix}\dot x^o\\\dot x^i\end{bmatrix}&=\begin{bmatrix}A&0\\\hat BH&\hat A\end{bmatrix}\begin{bmatrix}x^o\\x^i\end{bmatrix}+
\begin{bmatrix}B\\\hat BJ\end{bmatrix}u,\\
\begin{bmatrix}\dot x_1^o\\\dot x_2^o\\\dot x^i\end{bmatrix}&=\begin{bmatrix}A_{11}&A_{12}&0\\A_{21}&A_{22}&0\\\hat BH_1&\hat BH_2&\hat A\end{bmatrix}
\begin{bmatrix}x_1^o\\x_2^o\\x^i\end{bmatrix}+\begin{bmatrix}B_1\\B_2\\\hat BJ\end{bmatrix}u,\\
y&=\begin{bmatrix}UH&\hat C\end{bmatrix}\begin{bmatrix}x^o\\x^i\end{bmatrix}+UJu.
\end{align*}
Substituting from Green: \cite[(3.17)]{GreenInnerOuterSCL1988}, $A_{21}=-\hat BH_1$; \cite[(3.19)]{GreenInnerOuterSCL1988}, $B_2=-\hat BJ$; and the definition of $\hat A$ above
\begin{align*}
\begin{bmatrix}\dot x_1^o\\\dot x_2^o\\\dot x^i\end{bmatrix}%&=\begin{bmatrix}A_{11}&A_{12}&0\\A_{21}&A_{22}&0\\\hat BH_1&\hat BH_2&\hat A\end{bmatrix}
%\begin{bmatrix}x_{o1}\\x_{o2}\\x_i\end{bmatrix}+\begin{bmatrix}B_1\\B_2\\\hat BJ\end{bmatrix}u,\\
&=\begin{bmatrix}A_{11}&A_{12}&0\\A_{21}&A_{22}&0\\-A_{21}&\hat A-A_{22}&\hat A\end{bmatrix}
\begin{bmatrix}x_1^o\\x_2^o\\x^i\end{bmatrix}+\begin{bmatrix}B_1\\B_2\\-B_2\end{bmatrix}u,
\end{align*}
Applying the state transformation $\begin{bmatrix}x_1^o\\x_2^o\\x^i+x_2^o\end{bmatrix}=\begin{bmatrix}I&0&0\\0&I&0\\0&I&I\end{bmatrix}\begin{bmatrix}x_1^o\\x_2^o\\x_i\end{bmatrix}=\begin{bmatrix}I&0&0\\0&I&0\\0&-I&I\end{bmatrix}^{-1}\begin{bmatrix}x_1^o\\x_2^o\\x^i\end{bmatrix}$ yields
\begin{align}\label{eq:decomp}
\begin{bmatrix}\dot x_1^o\\\dot x_2^o\\\dot x^i+\dot x_2^o\end{bmatrix}&=\begin{bmatrix}A_{11}&A_{12}&0\\A_{21}&A_{22}&0\\0&0&\hat A\end{bmatrix}\begin{bmatrix}x_1^o\\x_2^o\\x^i+x_2^o\end{bmatrix}+\begin{bmatrix}B_1\\B_2\\0\end{bmatrix}u,
\end{align}
which establishes the result on unreachability of the $\ell$-dimensional subspace spanned by $x^i+x_2^o$ and the associated eigenvalues of $\hat A$.
\end{proof}

\section*{A.III State estimation variant of Green's algorithm: $P(s)=P_o(s)P_i(s).$}

For clarity, adopt the notation $P(s)=\mcbD+\mcbC(sI-\mcbA)^{-1}\mcbB$ for the state-variable realization of $P$. Then, applying Green's construction \cite{GreenInnerOuterSCL1988}  above to $P^\top$, 
\begin{align*}
P^\top(s)&=\mcbD^\top+\mcbB^\top(sI-\mcbA^\top)^{-1}\mcbC^\top,\\
&=P_i(s)P_o(s),\\
&=\left[U+\hat C(sI-\hat A)^{-1}\hat B\right].\left[J+H(sI-\mcbA^\top)^{-1}\mcbC^\top\right],\\
P(s)&=P^\top_o(s)P^\top_i(s),\\
&=\left[J^\top+\mcbC(sI-\mcbA)^{-1}H^\top\right].\left[U^\top+\hat B^\top(sI-\hat A^\top)^{-1}\hat C^\top\right].
\end{align*}
Write the state-variable realization of this latter product.
\begin{align}
\dot x_i&=\hat A^\top x_i+\hat C^\top u,\nonumber\\
f_i&=\hat B^\top x_i+U^\top u,\nonumber\\
\dot x_o&=\mcbA x_o+H^\top f_i,\\
&=\mcbA x_o+H^\top\hat B^\top x_i+H^\top U^\top u,\nonumber\\
y&=\mcbC x_o+J^\top f_i,\nonumber\\
&=\mcbC x_o+J^\top \hat B^\top x_i+J^\top U^\top u,\nonumber\\
\begin{bmatrix}\dot x_i\\\dot x_o\end{bmatrix}&=\begin{bmatrix}\hat A^\top&0\\H^\top\hat B^\top&\mcbA\end{bmatrix}\begin{bmatrix}x_i\\x_o\end{bmatrix}
+\begin{bmatrix}\hat C^\top\\H^\top U^\top\\\end{bmatrix}u,\label{eq:obstate}\\
y&=\begin{bmatrix}J^\top\hat B^\top&\mcbC\end{bmatrix}\begin{bmatrix}x_i\\x_o\end{bmatrix}+J^\top U^\top u.\label{eq:obout}
\end{align}

Here are Green's seven steps applied to $P^\top(s)$ with realization $(\mcbA^\top, \mcbC^\top,\mcbB^\top,\mcbD^\top)$.
\begin{enumerate}[label={[\arabic*]}]
\item \label{stp:1}Compute the observability Gramian, $Q$, of $P^\top$: $Q\mcbA^\top+\mcbA Q+\mcbB\mcbB^\top=0$.
\item \label{stp:2}Compute the spectral factor of $(P^\top)^\sim P^\top$: $P_o(s)=J+H(sI-\mcbA^\top)^{-1}\mcbC^\top.$
\item  \label{stp:3}Compute the observability Gramian, $X$, of $P_o$: $X\mcbA^\top+\mcbA X+H^\top H=0.$
\item  \label{stp:4}Find orthogonal $V$ so that: $V(Q-X)V^{-1}=\begin{bmatrix}0&0\\0&\Sigma\end{bmatrix}$.
\item  \label{stp:5}Transform and partition: $V\mcbA^\top V^{-1}=\begin{bmatrix}\mcbA_{11}^\top&\mcbA_{21}^\top\\\mcbA_{12}^\top&\mcbA_{22}^\top\end{bmatrix},$
$V\mcbC^\top=\begin{bmatrix}\mcbC_1^\top\\\mcbC_2^\top\end{bmatrix},$ $\mcbB^\top V^{-1}=\begin{bmatrix}\mcbB_1^\top&\mcbB_2^\top\end{bmatrix},$\\
$HV^{-1}=\begin{bmatrix}H_1&H_2\end{bmatrix}.$
\item  \label{stp:6}Find $U$ with $U^\top U=I$ and: $\begin{bmatrix}\mcbB_1^\top&\mcbD^\top\end{bmatrix}=U\begin{bmatrix}H_1&J\end{bmatrix}.$
\item  \label{stp:7}Inner factor: $\hat C=UH_2-\mcbB_2^\top,$ $\hat B=\Sigma^{-1}(\mcbB_2U-H_2^\top),$ $\hat A=\mcbA_{22}^\top+\hat BH_2$.
\end{enumerate}
From here, \eqref{eq:obstate}-\eqref{eq:obout} become
\begin{align*}
\begin{bmatrix}\dot x_i\\\dot x_o\end{bmatrix}&=\begin{bmatrix}\hat A^\top&0\\H^\top\hat B^\top&\mcbA\end{bmatrix}\begin{bmatrix}x_i\\x_o\end{bmatrix}
+\begin{bmatrix}\hat C^\top\\H^\top U^\top\\\end{bmatrix}u,\\
\begin{bmatrix}\dot x_i\\\dot x_{o1}\\\dot x_{o2}\end{bmatrix}&
=\begin{bmatrix}\hat A^\top&0&0\\H_1^\top\hat B^\top&\mcbA_{11}&\mcbA_{12}\\H_2^\top \hat B^\top&\mcbA_{21}&\mcbA_{22}\end{bmatrix}%\begin{bmatrix}x_i\\x_o\end{bmatrix}
%+\begin{bmatrix}\hat C^\top\\H^\top U^\top\end{bmatrix}u,\\
%&=\begin{bmatrix}\mcbA_{22}+H_2^\top(U^\top\mcbB_2^\top-H_2)\Sigma^{-1}&0\\H^\top(U^\top\mcbB_2^\top-H_2)\Sigma^{-1}&\mcbA\end{bmatrix}
\begin{bmatrix}x_i\\x_{o1}\\x_{o2}\end{bmatrix}
+\begin{bmatrix}\hat C^\top\\H_1^\top U^\top\\H_2^\top U^\top\end{bmatrix}u,\\
y&=\begin{bmatrix}J^\top\hat B^\top&\mcbC\end{bmatrix}\begin{bmatrix}x_i\\x_o\end{bmatrix}+J^\top U^\top u,\\
&=\begin{bmatrix}J^\top\hat B^\top&\mcbC_1&\mcbC_2\end{bmatrix}\begin{bmatrix}x_i\\x_{o1}\\x_{o2}\end{bmatrix}+J^\top U^\top u.
%&=\begin{bmatrix}J^\top(U^\top\mcbB_2^\top-H_2)\Sigma^{-1}&\mcbC\end{bmatrix}\begin{bmatrix}x_i\\x_o\end{bmatrix}+J^\top U^\top u.
\end{align*}
Make the following substitutions from the algorithm.\\
\centerline{
\begin{tabular}{|l|l||l|l|}
\hline
$H_1^\top\hat B^\top=-\mcbA_{12}$& \cite[(3.17)]{GreenInnerOuterSCL1988}&$H_1^\top U^\top=\mcbB_1$&Step~\ref{stp:6}\\
$H_2^\top\hat B^\top=\hat A^\top-\mcbA_{22}$&Step~\ref{stp:7}&$H_2^\top U^\top=\hat C^\top+\mcbB_2$&Step~\ref{stp:7}\\
$J^\top\hat B^\top=-\mcbC_2$&\cite[(3.19)]{GreenInnerOuterSCL1988}&$J^\top U^\top=\mcbD$&Step~\ref{stp:6}\\
\hline
\end{tabular}}
This yields
\begin{align*}
\begin{bmatrix}\dot x_i\\\dot x_{o1}\\\dot x_{o2}\end{bmatrix}&
=\begin{bmatrix}\hat A^\top&0&0\\-\mcbA_{12}&\mcbA_{11}&\mcbA_{12}\\\hat A^\top-\mcbA_{22}&\mcbA_{21}&\mcbA_{22}\end{bmatrix}%\begin{bmatrix}x_i\\x_o\end{bmatrix}
\begin{bmatrix}x_i\\x_{o1}\\x_{o2}\end{bmatrix}
+\begin{bmatrix}\hat C^\top\\\mcbB_1\\\hat C^\top+\mcbB_2\end{bmatrix}u,\\
y&=\begin{bmatrix}-\mcbC_2&\mcbC_1&\mcbC_2\end{bmatrix}\begin{bmatrix}x_i\\x_{o1}\\x_{o2}\end{bmatrix}+\mcbD u.
\end{align*}
Next transform the state
\begin{align*}
\begin{bmatrix}x_i\\x_{o1}\\x_{o2}-x_i\end{bmatrix}&=\begin{bmatrix}I&0&0\\0&I&0\\-I&0&I\end{bmatrix}\begin{bmatrix}x_i\\x_{o1}\\x_{o2}\end{bmatrix}
\end{align*}
Then
\begin{align}
\begin{bmatrix}\dot x_i\\\dot x_{o1}\\\dot x_{o2}-\dot x_i\end{bmatrix}&=\begin{bmatrix}I&0&0\\0&I&0\\-I&0&I\end{bmatrix}
\begin{bmatrix}\hat A^\top&0&0\\-\mcbA_{12}&\mcbA_{11}&\mcbA_{12}\\\hat A^\top-\mcbA_{22}&\mcbA_{21}&\mcbA_{22}\end{bmatrix}
\begin{bmatrix}I&0&0\\0&I&0\\I&0&I\end{bmatrix}\begin{bmatrix}x_i\\x_{o1}\\x_{o2}-x_i\end{bmatrix}+
\begin{bmatrix}I&0&0\\0&I&0\\-I&0&I\end{bmatrix}\begin{bmatrix}\hat C^\top\\\mcbB_1\\\hat C^\top+\mcbB_2\end{bmatrix}u,\nonumber\\
&=\begin{bmatrix}\hat A^\top&0&0\\0&\mcbA_{11}&\mcbA_{12}\\0&\mcbA_{21}&\mcbA_{22}\end{bmatrix}\begin{bmatrix}x_i\\x_{o1}\\x_{o2}-x_i\end{bmatrix}+
\begin{bmatrix}\hat C^\top\\\mcbB_1\\\mcbB_2\end{bmatrix}u,\label{eq:omgstate}\\
%&=\begin{bmatrix}\hat A^\top&0&0\\0&\mcbA_{11}&\mcbA_{12}\\0&\mcbA_{21}&\mcbA_{22}\end{bmatrix}\begin{bmatrix}x_i\\x_{o1}\\x_{o2}-x_i\end{bmatrix}+
%\begin{bmatrix}H_2^\top U^\top-\mcbB_2\\\mcbB_1\\\mcbB_2\end{bmatrix}u,\\
y&=\begin{bmatrix}-\mcbC_2&\mcbC_1&\mcbC_2\end{bmatrix}\begin{bmatrix}I&0&0\\0&I&0\\I&0&I\end{bmatrix}\begin{bmatrix}x_i\\x_{o1}\\x_{o2}-x_i\end{bmatrix}+J^\top U^\top u,\nonumber\\
&=\begin{bmatrix}0&\mcbC_1&\mcbC_2\end{bmatrix}\begin{bmatrix}x_i\\x_{o1}\\x_{o2}-x_i\end{bmatrix}+\mcbD u.\label{eq:omgout}
\end{align}

\section*{Remarks}
\begin{enumerate}[label=\roman*.]
\item The state of the inner factor, $x_i$, is unobservable. This means that the output measurements, $y_t$, cannot improve estimation quality versus the simulation-based estimator
\begin{align*}
\dot{\hat x}_i&=\hat A^\top\hat x_i+\hat C^\top u.
\end{align*}
See \cite{LiuBitmead_Autom2011} for an appreciation of stochastic observability.
\item In a stochastic formulation such as \eqref{eq:pxdef}-\eqref{eq:pydef} with process noise $w_t$ and given the analysis of \cite{BitmeadHovdAbooshahab:Autom2019} which demonstrates that SISE arises as the limit that $d_t$ is white noise of unbounded covariance, there is no benefit to: modeling $w_t$ as affecting the $x_i$ state component in a non-minimal description; or, in taking an initial $x_i$-state covariance as non-zero. Accordingly, $\hat x_i$ is exact.
\item For this non-minimal state-space realization \eqref{eq:omgstate}-\eqref{eq:omgout}, the state of the outer factor is completely observable and involves the estimate $\hat x_i$ computed as immediately above with $u_t=d_t$ and the measurements $y_t$ along with other known or modeled inputs signals.
\item From  \eqref{eq:omgstate}-\eqref{eq:omgout}, we recover the minimal state-space realization $P(s)=\mcbD+\mcbC(sI-\mcbA)^{-1}\mcbB$. Theorem~\ref{thm:same} shows that the state of this system is asymptotically  identical to the state of the outer factor.
\end{enumerate}

\section*{Acknowledgement} The author is pleased to acknowledge the valuable feedback from his collaborators Professor Morten Hovd, Norwegian University of Science \&\ Technology, and Dr Ali Abooshahab, Bouvet Norway.

\bibliographystyle{IEEEtran} %\bibliographystyle{ieeetr}        % Include this if you use bibtex 
\bibliography{/Users/bob/tex/bobilby} 
\end{document}